# THE FRATTINI SUBGROUP FOR SUBGROUPS OF HYPERBOLIC GROUPS

ILYA KAPOVICH

ABSTRACT. We prove that for a finitely generated subgroup $H$ of a word-hyperbolic group $G$ the Frattini subgroup $F(H)$ of $H$ is finite.

## 1. INTRODUCTION

The *Frattini subgroup* $F(G)$ of a group $G$ is defined as the intersection of all maximal subgroups of $G$, provided at least one maximal subgroup exists, and as $F(G) := G$ otherwise. It is easy to see that $F(G)$ is a characteristic subgroup of $G$ and hence normal in $G$.

Our main goal is to prove the following:

**Theorem A.** *Let $H$ be a finitely generated subgroup of a word-hyperbolic group $G$. Then the Frattini subgroup $F(H)$ of $H$ is finite.*

For a subset $S$ of a group $G$ we will denote by $\langle S \rangle$ the subgroup of $G$ generated by $S$. Recall that an element $g \in G$ is said to be a *non-generator* if for any generating set $S$ for $G$ we have $\langle S - \{g\}\rangle = G$. We shall need the following well-known alternative characterization of the Frattini subgroup (see for example [25]).

**Lemma 1.1.** *Let $G$ be a group. Then*
$$F(G) := \{g \in G \,|\, g \text{ is a non-generator of } G\}.$$

Our interest in the Frattini subgroup for hyperbolic groups is primarily motivated by the connections with the *generation problem* and the *rank problem*. The generation problem asks for a finite subset $Y$ of a fixed group $G$ whether $Y$ generates $G$. The rank problem for a class of finitely generated groups asks, given a group $G$ in the class, what is the smallest cardinality of a generating set for $G$ (called the *rank* of $G$ and denoted $rk(G)$). Both the rank problem and the generation problem are known to be unsolvable for the classes of hyperbolic and torsion-free hyperbolic groups [8], but solvable for the class of torsion-free locally quasiconvex hyperbolic groups [14, 20]. Note that finitely generated subgroups of hyperbolic groups need not be finitely presentable [24] and that finitely presentable subgroups of hyperbolic groups need not be hyperbolic [5].

Having a large Frattini subgroup often allows one to reduce the rank and the generation problems to some simpler quotient groups. Thus if $G$ and $F(G)$ are finitely generated, then by Lemma 1.1:

1. For a subset $Y \subseteq G$ we have $\langle Y \rangle = G \iff \langle \overline{Y} \rangle = G/F(G)$, where $\overline{Y}$ is the image of $Y$ in the quotient group $G/F(G)$.
2. $rk(G) = rk(G/F(G))$.

Thus if $G$ is a finitely generated nilpotent group then $F(G)$ contains the commutator subgroup of $G$ (see for example [25]). This implies that a subset $Y \subseteq G$ generates $G$ if and only if it generates $G$ modulo $[G, G]$. Hence the rank of $G$ is equal to the rank of the abelianization of $G$ and so $rk(G)$ is computable.

Theorem A essentially shows that no such help is forthcoming for any reasonable class of hyperbolic groups. Moreover, a residually finite group obviously always has a wealth of maximal subgroups of finite index. Although it is unknown if all hyperbolic groups are residually finite, the conventional wisdom asserts that there exists a non-residually finite hyperbolic group. This would imply by a result of I.Kapovich and D.Wise [21] that there exists an infinite hyperbolic group $G$ which has no proper subgroups of finite index. Yet according to Theorem A the group $G$ still has a rich collection of maximal subgroups whose intersection is finite (and even trivial if $G$ is torsion-free). This indicates that most maximal subgroups in hyperbolic groups do not corresponds to pull-backs of maximal subgroups from finite quotients and thus have rather pathological nature. A typical example is provided by the result of A.Ol'shanskii [23] which asserts that every non-elementary torsion-free hyperbolic group $G$ possesses an infinite non-abelian torsion-free quotient $\bar{G}$ such that every proper subgroup of $\bar{G}$ is cyclic. Then the pull-back of any maximal cyclic subgroup of $\bar{G}$ is a maximal subgroup of $G$. We should also note that by a well-known result of B.Wehrfritz [28] finitely generated linear groups are known to have nilpotent Frattini subgroups. A nilpotent subgroup of a hyperbolic group is necessarily virtually cyclic. Moreover, if a virtually cyclic subgroup of a hyperbolic group is normal then this subgroup is either finite or has finite index. Since it is easy to see that virtually cyclic groups have finite Frattini subgroups (see Lemma 2.1 below), it follows that linear hyperbolic groups have finite nilpotent Frattini subgroups. Apart from this fact, it seems that little had been previously known about Frattini subgroups of hyperbolic groups (see [4, 2] for some results applicable to hyperbolic 3-manifold and surface groups).

It is worth pointing out that Theorem A implies that the Frattini subgroup is trivial for a much larger class of groups than torsion-free hyperbolic groups and their finitely generated subgroups. Recall that if $\mathcal{C}$ is a class of groups (not necessarily closed under taking subgroups) then a group $G$ is said to be *residually-$\mathcal{C}$* if for any $g \in G, g \neq 1$ there exists a homomorphism $\phi : G \to C$, where $C \in \mathcal{C}$ and $\phi(g) \neq 1$.

**Corollary B.** *Let $G$ be a finitely generated group which is residually torsion-free word-hyperbolic (e.g. residually free). Then $F(G) = 1$.*

*Proof.* It is easy to see from the definition that if $\phi : G_1 \to G_2$ is an epimorphism, then $F(G_1) \leq \phi^{-1}(F(G_2))$ since the full pre-image of a maximal subgroup of $G_2$ is a maximal subgroup of $G_1$. Since for finitely generated torsion-free subgroups of hyperbolic groups the Frattini subgroup is trivial by Theorem A, this immediately implies the statement of Corollary B. □

Residually hyperbolic and, in particular, residually free groups have been the object of intensive study in recent years [26, 17, 18, 6, 7].

The author is grateful to Derek Robinson, Paul Schupp, Peter Brinkmann and Bogdan Petrenko for useful discussions and to Sergei Ivanov and Brad Edge for providing the inspiration to consider the problem addressed by this article.

## 2. Preliminary facts about hyperbolic groups

Recall that every subgroup $H$ of a hyperbolic group $G$ is either virtually cyclic (in which case $H$ is called *elementary*) or contains a free non-abelian group of rank two (in which case $H$ is called *non-elementary*).

**Lemma 2.1.** *Let $H$ be an elementary subgroup of a hyperbolic group $G$. Then $F(H)$ is finite.*

*Proof.* If $H$ is finite then $F(H)$ is clearly finite. Thus we may assume that $H$ is infinite and hence virtually infinite cyclic. Therefore there exists an epimorphism $\phi : H \to C$ such that $N = ker(\phi)$ is finite and such that $C$ is either infinite cyclic or infinite dihedral. It is easy to see that both the infinite cyclic and infinite dihedral groups have trivial Frattini subgroups. Hence $F(H) \leq \phi^{-1}(F(C)) = \phi^{-1}(1) = N$ is finite, as required. □

For the remainder of this article, unless specified otherwise, let $G$ be a non-elementary word-hyperbolic group with a fixed finite generating set $A$ and the Cayley graph $X := \Gamma(G, A)$. Let $d$ be the word-metric on $X$ corresponding to $A$. Let $\delta \geq 0$ be an integer such that $X$ is $\delta$-*hyperbolic*, that is for any geodesic triangle in $X$ each side of this triangle is contained in the closed $\delta$-neighborhood of the union of the other sides. For $g \in G$ we will denote $|g|_A := d(1, g)$. We shall often denote a geodesic segment from $x \in X$ to $y \in X$ by $[x, y]$. Recall that a subgroup $H \leq G$ is said to be *quasiconvex* if there exists $E > 0$ such that for any $h_1, h_2 \in H$ any geodesic segment $[h_1, h_2]$ in $X$ is contained in the closed $E$-neighborhood of $H$. We refer the reader to [12, 3, 9, 13, 19, 22, 27] for the background information on word-hyperbolic groups and their quasiconvex subgroups.

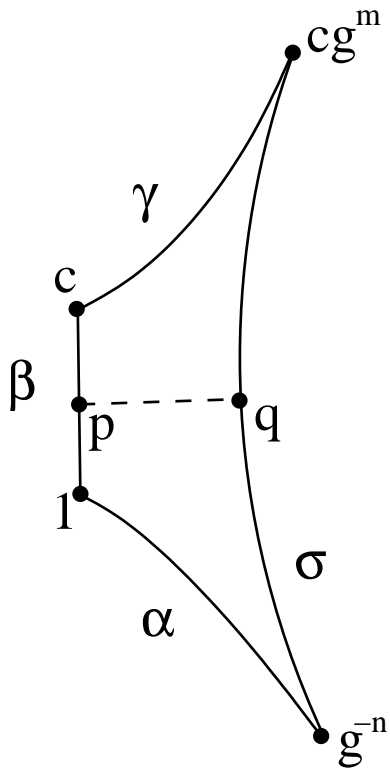

FIGURE 1. The case $d(\beta, \sigma) \leq 2\delta$

We recall the following known facts regarding infinite cyclic subgroups and their commensurators in hyperbolic groups (see for example [23, 19]):

**Proposition 2.2.** *Let $g \in G$ be an element of infinite order. Then*

1. *Suppose $h^{-1}g^n h = g^m$ for some integers $m, n$ and some $h \in G$. Then $|m| = |n|$.*
2. *The set of elements*

$$E_G(g) := E(g) := \{h \in G \mid \text{ for some } n \neq 0 \quad h^{-1}g^n h = g^n \text{ or } h^{-1}g^n h = g^{-n}\}$$

   *forms a subgroup of $G$ containing $\langle g \rangle$ as a subgroup of finite index.*
3. *The subgroup $E(g)$ is equal to the commensurator of $\langle g \rangle$ in $G$, that is*

$$E(g) = \{h \in G \mid [\langle g \rangle : \langle g \rangle \cap h^{-1}\langle g \rangle h] < \infty, \ [h^{-1}\langle g \rangle h : \langle g \rangle \cap h^{-1}\langle g \rangle h] < \infty\}.$$

4. *The subgroups $\langle g \rangle$ and $E(g)$ are quasiconvex in $G$.*

We shall reserve the notation $E(g) = E_G(g)$ for the commensurator of the cyclic subgroup $\langle g \rangle$ in $G$. We recall the following useful lemma due to T.Delzant (see Lemma 1.1 in [10]):

**Lemma 2.3.** *Let $a > 0$ and suppose that $(x_n)_{n \in J}$ is a sequence of points in $X$ (where $J$ is a subinterval of $\mathbb{Z}$ consisting of at least three numbers) such that*

$$d(x_{n+2}, x_n) \geq \max\{d(x_{n+2}, x_{n+1}), d(x_{n+1}, x_n)\} + 2\delta + a$$

*whenever $n, n+2 \in J$.*

*Then $d(x_n, x_p) \geq a|n - p|$ whenever $n, p \in J$.*

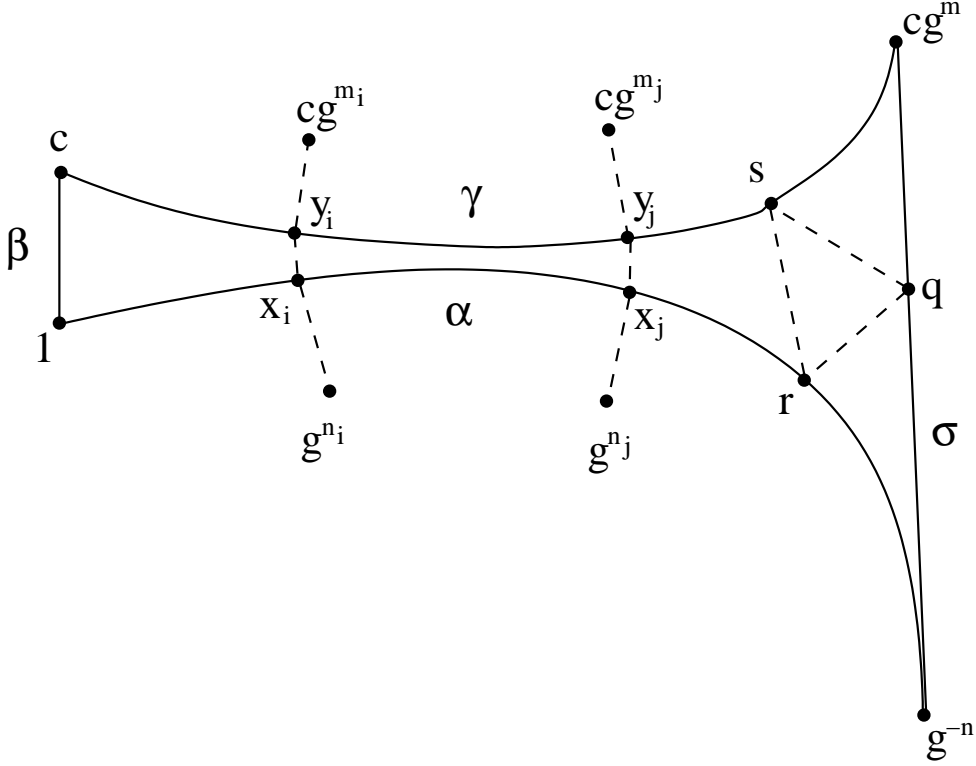

FIGURE 2. The case $d(\beta, \sigma) > 2\delta$

**Proposition 2.4.** *Let $g \in G$ be an element of infinite order and suppose $c \in G$ is such that $c \notin E(g)$. Then there exists a constant $K = K(g, c) > 0$ such that for any integers $m, n$*

$$|g^n c g^m|_A \geq |g^n|_A + |g^m|_A - K.$$

*Proof.* The proof is a fairly standard hyperbolic exercise and is similar to the arguments used in [11, 1, 15, 16]. We present the details for completeness.

Let $E > 0$ be the quasiconvexity constant of the subgroup $H := \langle g \rangle$ in $G$. Consider a geodesic quadrilateral $\Delta$ in $X$ with vertices $1, c, cg^m, g^{-n}$ and geodesic sides $\alpha = [1, g^{-n}]$, $\beta = [1, c]$, $\gamma = c[1, g^m]$ and $\sigma = g^{-n}[1, g^n c g^m] = [g^{-n}, cg^m]$.

Note that each side of $\Delta$ is contained in the closed $2\delta$-neighborhood of the three other sides since $X$ is $\delta$-hyperbolic. Suppose first that $d(\beta, \sigma) \leq 2\delta$, as shown in Figure 2. Let $p \in \beta$ and $q \in \sigma$ be such that $d(p, q) \leq \delta$. Then by the triangle inequality $d(g^{-n}, q) \geq d(1, g^{-n}) - d(1, p) \geq |g^n|_A - |c|_A$. Similarly, $d(q, cg^m) \geq |g^m|_A - |c|_A$. Hence

$$|g^n c g^m|_A = d(g^{-n}, cg^m) = d(g^{-n}, q) + d(q, cg^m) \geq |g^n|_A + |g^m|_A - 2|c|_A.$$

Suppose now that $d(\beta, \sigma) > 2\delta$. Then every point of $\sigma$ is contained in the closed $2\delta$-neighborhood of either $\alpha$ or $\gamma$. Hence, since $\sigma$ is connected, there is a point $q \in \sigma$ such that $d(q, \alpha) \leq 2\delta$ and $d(q, \gamma) \leq 2\delta$, as shown in Figure 2. Let $s \in \gamma$ and $r \in \alpha$ be such that $d(q, s) \leq 2\delta$ and $d(q, r) \leq 2\delta$ and therefore $d(r, s) \leq 4\delta$. Let $N$ be the number of elements in $G$ of length at most $2E + 2|c|_A + 4\delta$.

**Claim.** We have $d(1, r) < (2E + 1)N$.

Indeed, suppose that $d(1, r) \geq (2E+1)N$. Since $d(q, r) \leq 4\delta$, for any point $x \in \alpha$ with $d(1, x) \leq d(1, r)$ there is a point $y \in \gamma$ with $d(c, y) \leq d(c, s)$ such that $d(x, y) \leq |c|_A + 4\delta$.

Consider a sequence of points $x_0 = 1, x_1 \ldots, x_N$ on $\alpha$ so that $d(1, x_i) = i(2E + 1)$. Since $H$ is $E$-quasiconvex in $X$, for each $i = 1, \ldots, N$ there is $n_i$ such that $d(x_i, g^{n_i}) \leq E$. Put $n_0 = 0$, so that $d(x_0, g^{n_0}) = d(1, 1) = 0$. Since by assumption $d(1, r) \geq N(2E + 1)$, all points $x_i$ belong to the segment of $\alpha$ between 1 and $r$. For each $i = 0, \ldots, N$ there is a point $y_i$ on $\gamma$ such that $d(x_i, y_i) \leq |c|_A + 4\delta$ and that $d(c, y_i) \leq d(c, s)$. Again, since $H$ is $E$-quasiconvex, for each $i = 0, \ldots, N$ there is an integer $m_i$ such that $d(y_i, cg^{m_i}) \leq E$. Hence $d(g^{n_i}, cg^{m_i}) \leq 2E + |c|_A + 4\delta$. By the choice of $N$ this means that there exist $i < j$ such that

$$g^{-n_i} cg^{m_i} = g^{-n_j} cg^{m_j} = f \in G.$$

Hence $f^{-1} g^{n_j - n_i} f = g^{m_j - m_i}$. Note that by construction $d(g^{n_i}, g^{n_j}) \geq d(x_i, x_j) - 2E = (j - i)(2E + 1) - 2E > 0$. Hence $n_j - n_i \neq 0$ and therefore by Proposition 2.2 $f \in E_G(g)$. Hence $c = g^{n_i} f g^{-m_i} \in E_G(g)$, contrary to our assumptions.

Thus the Claim is proved and $d(1, r) < (2E + 1)N$. Since $d(r, s) \leq 4\delta$, this implies that $d(x, s) \leq (2E + 1)N + |c|_A + 4\delta$. Then

$$d(cg^m, q) \geq d(cg^m, s) - 2\delta \geq |g^m|_A - (2E + 1)N - |c|_A - 6\delta \text{ and}$$
$$d(g^{-n}, q) \geq d(g^{-n}, r) - 2\delta \geq |g^n|_A - (2E + 1)N - 2\delta.$$

Hence

$$|g^n c g^m|_A = d(g^{-n}, cg^m) = d(g^{-n}, q) + d(cg^m, q) \geq |g^n|_A + |g^m|_A - 2N(2E + 1) - |c|_A - 8\delta.$$

Thus the statement Proposition 2.4 holds with $K(g, c) := 2N(2E + 1) + 2|c|_A + 8\delta$. □

*Proof of Theorem A.* Let $G$ be a word-hyperbolic group and let $H \leq G$ be a finitely generated subgroup. By Lemma 2.1 we may assume that $H$ is non-elementary. Fix a finite generating set $A$ of $G$, the Cayley graph $X := \Gamma(G, A)$ and the corresponding word-metric $d$ on $X$. Also fix an integer $\delta \geq 0$ such that $X$ is $\delta$-hyperbolic. Since any torsion subgroup of a hyperbolic group is finite [13, 9, 23], it suffices to show that $F(H)$ is a torsion group.

Let $g \in H$ be an arbitrary element of infinite order. It is enough to prove that $g \notin F(H)$, in other words that $g$ is not a non-generator for $H$.

Let $S = \{s_1, \ldots, s_t\}$ be a generating set for $H$. Since $H$ is non-elementary and $H$ is not contained in $E_G(g)$, there is some $i$ such that $s_i \notin E_G(g)$. We will assume that $s_1 \notin E_G(g)$. For those $2 \leq j \leq k$ with $s_j \in E_G(g)$ put $c_j := s_j s_1$. For those $1 \leq j \leq t$ with $s_j \notin E_G(g)$ put $c_j := s_j$. Then the set $S' = \{c_1, \ldots, c_t\}$ generates $H$. Moreover, for $j = 1, \ldots, t$ we have $c_j \notin E_G(g)$. Let

$$T_0 := \min\{|f|_A \text{ where } f \in H, f \neq 1\}$$

and put $T := T_0 + 1$.

Let $K_1 := \max\{K(g, c_j) \mid j = 1, \ldots, t\}$, where $K(g, c_j)$ is the constant provided by Proposition 2.4. Put $K_2 := \max\{|c_j|_A \text{ where } j = 1, \ldots, t\}$. Finally, put $K := \max(K_1, K_2)$. Since $\langle g \rangle \leq G$ is quasiconvex in $G$, there is $C > 0$ such that for any $n \in \mathbb{Z}$ we have $d(1, g^n) = |g^n|_A \geq C|n|$. Let $N > 1$ be an integer such that $CN \geq 3K + 2\delta + 100T$. Put $h_j := g^{n_j} c_j g^{10n_j}$, $j = 1, \ldots, t$, where $n_j := 1000jN$ for $j = 1, \ldots, t$. Then the set $Q := \{g, h_1, \ldots, h_t\}$ generates $H$. We will show that $Q' := Q - \{g\} = \{h_1, \ldots, h_t\}$ does not generate $H$.

Indeed, suppose $w = y_1 \ldots y_q$ is an arbitrary nontrivial freely reduced product in $Q'$, so that $q \geq 1$ and $y_i = h_{j_i}^{\epsilon_i}$, where $\epsilon_i = \pm 1$ and whenever $j_{i+1} = j_i$, then $\epsilon_{i+1} \neq -\epsilon_i$. Thus $y_j = g^{m_j} b_j g^{l_j}$ where $b_j = c_j^{\epsilon_i}$ and where $m_j = n_j, l_j = 10n_j$ if $\epsilon_j = 1$ and $m_j = -10n_j, l_j = -n_j$ if $\epsilon_j = -1$. Recall that $|b_j|_A \leq K$.

**Claim.** $|w|_A \geq Tq$.

If $q = 1$ and $w = y_1 = h_j^{\epsilon} = (g^{n_j} c_j g^{10n_j})^{\epsilon_1}$ then by Proposition 2.4 and the choice of $N$ we have

$$|w|_A \geq |g^{Nj}|_A + |g^{10Nj}|_A - K \geq CN - K \geq T,$$

as required. Suppose now that $q \geq 2$. Then we have the following equality in $G$:

$$w = g^{m_1} b_1 g^{l_1 + m_2} b_2 g^{l_2 + m_3} \ldots g^{l_{q-1} + m_q} b_q g^{l_q}.$$

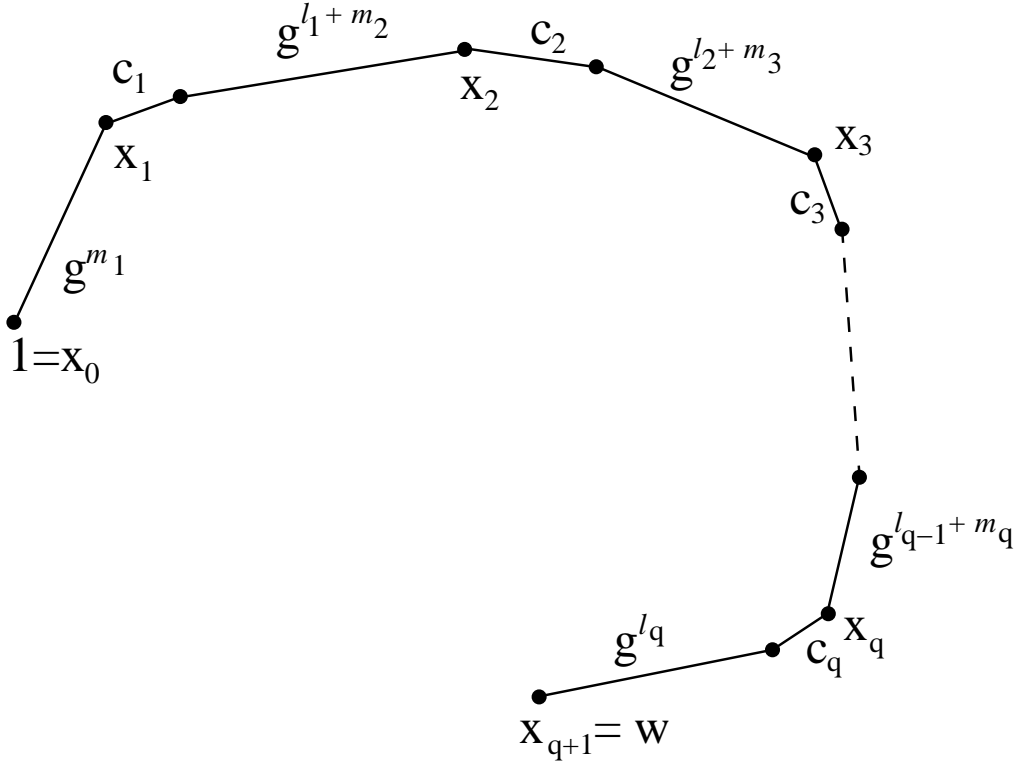

FIGURE 3. Application of Delzant's lemma

By the definition of $h_1, \ldots, h_t$ and the choice of $w$ we have $|l_i + m_{i+1}| \geq N$ and $|l_i| \geq N, |m_i| \geq N$. We consider the following sequence of points in $X$, as shown in Fifure 3:

$$x_0 := 1, x_1 := g^{m_1}, x_i = g^{m_1} b_1 g^{l_1+m_2} b_2 \ldots g^{l_{i-1}+m_i} \text{ for } i = 1, \ldots, q, \quad \text{and } x_{q+1} := w.$$

We will show that Lemma 2.3 applies to the sequence $x_0, x_1, \ldots, x_{q+1}$ in $X$ with $a = T$. Indeed, for $x_0 = 1, x_1 = g^{m_1}, x_2 = g^{m_1} b_1 g^{l_1+m_2}$ we have

$$d(x_0, x_1) = |g^{m_1}|_A, \quad |g^{l_1+m_2}|_A + K \geq d(x_1, x_2) = |b_1 g^{l_1+m_2}|_A \geq |g^{l_1+m_2}|_A - K.$$

Also, by Proposition 2.4 we have

$$d(x_0, x_2) = |g^{m_1} b_1 g^{l_1+m_2}| \geq |g^{m_1}|_A + |g^{l_1+m_2}|_A - K \geq \max\{|g^{m_1}|_A, |g^{l_1+m_2}|_A\} + CN - K \geq$$
$$\geq \max\{d(x_0, x_1), d(x_1, x_2)\} + CN - 2K \geq \max\{d(x_0, x_1), d(x_1, x_2)\} + T + 2\delta$$

since $|m_1|, |l_1 + m_2| \geq N$ and by the choice of $N$ we have $CN > 2K + 2\delta + T$. A similar argument shows that the conditions of Lemma 2.3 with $a = T$ hold for $x_{q-1}, x_q, x_{q+1}$.

Suppose now that $1 \leq i \leq q - 2$ and consider the points $x_i, x_{i+1}, x_{i+2}$. By construction $x_{i+1} = x_i b_i g^{l_i+m_{i+1}}$ and $x_{i+2} = x_i b_i g^{l_i+m_{i+1}} b_{i+1} g^{l_{i+1}+m_{i+2}}$. Hence

$$|g^{l_i+m_{i+1}}|_A + K \geq d(x_i, x_{i+1}) \geq |g^{l_i+m_{i+1}}|_A - K, \text{ and}$$
$$|g^{l_{i+1}+m_{i+2}}|_A + K \geq d(x_{i+1}, x_{i+2}) \geq |g^{l_{i+1}+m_{i+2}}|_A - K.$$

Moreover, by Proposition 2.4, we have

$$d(x_i, x_{i+2}) \geq |g^{l_i+m_{i+1}}|_A + |g^{l_{i+1}+m_{i+2}}|_A - 2K \geq \max\{|g^{l_i+m_{i+1}}|_A, |g^{l_{i+1}+m_{i+2}}|_A\} + CN - 2K \geq$$
$$\max\{d(x_i, x_{i+1}), d(x_{i+1}, x_{i+2})\} + CN - 3K \geq \max\{d(x_i, x_{i+1}), d(x_{i+1}, x_{i+2})\} + T + 2\delta,$$

as required. Therefore by Lemma 2.3 $|w|_A = d(x_0, x_{q+1}) \geq T(q+1) \geq Tq$ and the Claim is established.

Hence we have seen that for any nontrivial element $h$ from the subgroup $H_1 := \langle Q - \{g\}\rangle = \langle h_1, \ldots, h_t\rangle$ we have $|h|_A \geq T > T_0$. Hence $H_1 \neq H$ (since $H$ has an element of length $T_0$) and thus $g$ is not a non-generator for $H$. Therefore by Lemma 1.1 $g \notin F(H)$, as required. Since $g$ was chosen as an arbitrary element of infinite order in $H$, this implies that $F(H)$ is a torsion group and hence is finite. This completes the proof of Theorem A. □

Department of Mathematics, University of Illinois at Urbana-Champaign, 1409 West Green Street, Urbana, IL 61801, USA
*E-mail address*: `kapovich@math.uiuc.edu`